\documentclass[12pt]{amsart}               
\usepackage{amsmath,amssymb}
\usepackage{amscd}      
\usepackage{mathrsfs}     
\usepackage[all]{xy}
\usepackage{pstcol}
\UseComputerModernTips

\usepackage{graphicx} 
\DeclareGraphicsExtensions{.eps} {\bf }

\theoremstyle{plain}          
\newtheorem{theorem}{Theorem}[subsection]

\newtheorem{lemma}{Lemma}[subsection]

\newtheorem{corollary}{Corollary}[subsection]
\newtheorem*{main}{Main~Theorem}
\newtheorem{proposition}{Proposition}[subsection]

\newcommand{\edge}[1]{\ar@{-}[#1]}

\theoremstyle{definition}
\newtheorem{definition}{Definition}[subsection]
\newtheorem{example}{Example}[subsection]

\theoremstyle{remark} 
\newtheorem{remark}{Remark}[subsection]

\def\dim{\operatorname{dim}}

\def\ker{\operatorname{ker}}
\def\Def{\operatorname{Def}}               
\def\Res{\operatorname{Res}}
\def\PRes{\operatorname{PRes}}

\def\coker{\operatorname{coker}}     

\def\Tr{\operatorname{Tr}}  
\def\im{\operatorname{im}} 
\def\cC{\mathbf C} 
\def\C{\mathbb C}

\def\E{\mathcal    E}   

\def\P{\mathbb     P}    
\def\Q{\mathbb    Q}     

\numberwithin{equation}{section}
\begin{document}
\title[finite representations of  a quiver]{ Finite representations of
a   quiver   arising   from   string   theory}   
\author{Xinyun   Zhu}
\address{Department of Mathematics\\ Oklahoma State University\\ Stillwater, OK 74078}
\address{Department of Mathematics\\ University of Illinois at Urbana-
Champaign\\ Urbana, IL 61801} \email{zhu5@math.uiuc.edu} \date{\today}

\begin{abstract}  Inspired by Cachazo, Katz and Vafa (``Geometric
transitions and $\mathcal  {N}=1$ quiver theories'' (hep-th/0108120)),
we examine representations of  ``${N}=1$ quivers'' arising from string
theory.  We derive some  mathematical consequences of the physics, and
show  that these  results are  a  natural extension  of Gabriel's  ADE
theorem.    Extending  the   usual  ADE   case  that   relates  quiver
representations  to curves  on surfaces,  we relate  these  new quiver
representations to curves on threefolds.
\end{abstract}
\maketitle

\section{Introduction}
In 1972, Gabriel \cite{GA} published his celebrated theorem about the
finite representation type of ADE quivers without relations.  Since
then the study of finite representation type quiver theory has been an
important topic because it provides a successful way to solve problems
in the representation theory of algebras and Lie groups.  Recently, it
has attracted the attention of physicists, (see \cite{CKV}) due to its
close relation with the study of D-branes.  A special type of quiver
arising from string theory, which we will call an ``$N=1$ ADE
quiver'', was introduced in \cite{CKV}(see
Definition~\ref{N=1-A-D-E-Quiver}).  This $N=1$ ADE quiver has a close
relation with the usual ADE quiver. The representations of $N=1$ ADE
quivers will satisfy the following important relations
(\ref{relation}),
\begin{equation*}
\sum_is_{ij}Q_{ji}Q_{ij}+P'_j(\Phi_j)=0,                          \quad
Q_{ij}\Phi_j=\Phi_iQ_{ij}=0.
\end{equation*}

The purpose of this paper is to construct, under certain conditions, a
finite-to-one correspondence between  the simple representations of an
$N=1$ ADE quiver and the positive  roots of the usual ADE quiver; this
matches the physicists' predictions.

The ``reflection functors'' which were used in \cite{BGP} to reprove
Gabriel's Theorem provide us a way to attack this problem.  In this
paper, we first modify the ``reflection functors'' in \cite{BGP}, and
then apply our modified reflectors to get our Main~Theorem in Section~
\ref{applyref}.  Some related results using different methods were given in 
\cite{Szendroi}.
 
This paper is organized as follows.  In Section~ \ref{section1}, we
briefly review the the geometry of threefolds with small resolutions.
In Section~ \ref{section2}, we give the definition of $N=1$ ADE
quivers and their representations, and we introduce our modified
reflection functors.  In Section~\ref{section3}, we apply our modified
reflection functors to get a proof of the Main~Theorem in
Section~\ref{applyref}. In Section~\ref{section5}, I give a
correspondence between indecomposable representations and ADE
configuration of curves.

{\bf Acknowledgements.} The contents of this article derives from a part of
my doctoral dissertation. I wish to thank my thesis adviser, Professor
Sheldon Katz, for his help and encouragement. It was an honor and a
privilege to have been his student, and I will be forever grateful for
his mathematical assistance and inspiration, moral and financial
support, unfailing kindness and inexhaustible patience.


\section{ A quick review of the geometry of  threefolds  for a general ADE
fibration over  $\cC.$}\label{section1} In this section,  we refer the
reader to \cite{CKV} and \cite{KM92}.   Let $C\subset Y$ be a rational
curve (not necessarily irreducible) in a 3-fold $Y$ with $K_Y$ trivial
in a neighborhood  of $C$ and $\pi:Y\to X$  a birational morphism such
that $\pi(C)=p\in  X$ and $\pi|_{Y-C}$  is an isomorphism  onto $X-p$.
We consider  an analytic neighborhood  of $p$ (still denoted  $X$) and
its inverse image under $\pi$ (still denoted $Y$).  By a lemma of Reid
\cite{Re83} (1.1, 1.14), the generic hyperplane section through $p$ is
a surface $X_0$ with an isolated rational double point, and the proper
transform of  $X_0$ is  a partial resolution  $Y_0\to X_0$  (i.e.  the
minimal resolution $Z_0\to X_0$ factors through $Y_0\to X_0$).

The  partial  resolution $Y_0\to  X_0$  determines combinatorial  data
$\Gamma_0\subset \Gamma$ consisting of  an ADE Dynkin diagram $\Gamma$
(the type of the singularity  $p$) and a subgraph $\Gamma_0$ (the dual
graph of the exceptional set of $Y_0$).

Let  $\mathcal   Z\to  \Def(Z_0),\,  \mathcal   Y\to  \Def(Y_0),$  and
$\mathcal  X\to \Def(X_0)$ be  semi-universal deformations  of $Z_0,\,
Y_0,$ and $X_0.$ Following \cite{KM92}, there are identifications
\begin{eqnarray}\label{deformation}
&&  \Def(Z_0)\cong V=:  \Res(\Gamma)\\ &&  \Def(Y_0)\cong V/{\mathfrak
W}_0=:  \PRes(\Gamma, \Gamma_0)\\ &&  \Def(X_0)\cong V/\mathfrak  W =:
\Def(\Gamma)
\end{eqnarray}
where  $V$  is the  complex  root  space  associated to  $\Gamma$  and
$\mathfrak W$ is its  Weyl group. ${\mathfrak W}_0\subset \mathfrak W$
is  the  subgroup  generated   by  reflections  of  the  simple  roots
corresponding to $\Gamma-\Gamma_0.$ Deformations of $Z_0$ or $Y_0$ can
be blown down to give  deformations of $X_0$ (\cite{Wahl} Theorem 1.4)
and the induced  classifying maps are given by  the natural maps $V\to
V/\mathfrak  W$ and  $V/{\mathfrak  W}_0\to V/\mathfrak  W$ under  the
above identifications.

  We can  view $X$ as  the total space  of a 1-parameter  family $X_t$
defined by the classifying map
\[g:\Delta\to \Def (\Gamma).\]
where $\Delta\subset \C$ is a small disk.
Similarly, we get the compatible family $Y_t$ given by a map
\[f:\Delta\to \PRes (\Gamma, \Gamma_0).\]

That is, we get the diagram
\[\xymatrix{\mathcal Z \ar[r]^{\tilde \sigma}\ar[d] & \mathcal
Y\ar[r]^{\tilde         \rho}\ar[d]&         \mathcal        X\ar[d]\\
\Res(\Gamma)\ar[r]^\sigma    &   \PRes(\Gamma,\Gamma_0)\ar[r]^\rho   &
\Def(\Gamma)\\ & \Delta\ar[u]^f \ar[ur]^g\\ }
\]
By   \cite{KM92},   $\mathcal   Y$    is   a   blowup   of   $\mathcal
{X}\times_{V/{\mathfrak W}}  V/{\mathfrak W}_0$ and $\mathcal  Z$ is a
blowup of $\mathcal {X}\times_{V/  {\mathfrak W}} V.$ By construction,
$Y$ is the pullback of $\mathcal Y$  by $f$ and $X$ is the pullback of
$\mathcal X$ by $g.$

\subsection{The geometry of  threefolds  with small resolutions  for a
general  ADE  fibration over  $\C$}  \label{fibration} Let  $X\subset
\C\times {\C}^3$ be  an ADE fibration over $\C$.   Let $t_i: \C\to
V$ be a map ($1\le i\le n+1$ in $A_n$ case, $1\le i\le n$ in $D_n$ and
$E_n$  case),  where  $V$  is   the  complex  root  space  defined  on
(\ref{deformation}). We consider the $A_n$ case first.  Then $X\subset
\C\times \C^3$ is defined by the equation
\[xy=z^{n+1}+\alpha_2(\omega)z^{n-1}+\cdots+\alpha_{n+1}(\omega).\]We
write $h: \C\to V\subset \C^{n+1}$ as
\begin{equation*}\label{CtoV}
h=(t_1, \cdots, t_{n+1}):\C\to \C^{n+1},\, \sum_{i=1}^{n+1}t_i=0.
\end{equation*}
Referring  to   \cite{KM92},  $\alpha_1,  \cdots,   \alpha_{n+1}$  are
elementary symmetric functions in $t_i, \cdots, t_{n+1}.$

Let $Z$ be the closure of the graph of the rational map
\[X\to (\P^1)^n, \, (x,y,z,\omega)\to
\left[x,\prod^i_{j=1}(z+t_{j}(\omega))     \right]_i.\]     and    let
$(u_i,v_i)$ be  coordinates of the $i$-th $\P^1$  in $(\P^1)^n.$ Using
the identities
$$\left      [x,      z+t_1(\omega)\right]=\left[(z+t_2(\omega))\cdots
(z+t_n(\omega)), -y\right],$$ we get
\[xv_j=u_j\prod^j_{i=1}(z+t_i(\omega))\, (1\le j\le n),\]
and
\[\prod^j_{i=k+1}(z+t_{i}(\omega))u_jv_k=u_kv_j \, (1\le k< j\le n).\]
 We  refer  the  reader   to  \cite{KM92}  for  the  more  complicated
fibrations of the $D$ and $E$  cases. We list the defining equation of
their deformations as follows:
\[D_n:\, x^2+y^2z+\frac{\prod^n_{i=1}(z+t^2_i(\omega))-\prod^n_{i=1}t^2_i
(\omega)}{z}+2\prod^n_{i=1}t_i(\omega)y\]
\[E_6:\,
x^2+z^4+y^3+\epsilon_2yz^2+\epsilon_5yz+\epsilon_6z^2+\epsilon_8y+
\epsilon_9z+\epsilon_{12}\]
\[E_7:\,
x^2+y^3+yz^3+\epsilon_2y^2z+\epsilon_6y^2+\epsilon_8yz+\epsilon_{10}
z^2+\epsilon_{12}y+\epsilon_{14}z+\epsilon_{18}\]
\[E_8:\,
x^2+y^3+z^5+\epsilon_2yz^3+\epsilon_8yz^2+\epsilon_{12}z^3+\epsilon_
{14}yz+\epsilon_{18}z^2+\epsilon_{20}y+\epsilon_{24}z+\epsilon_{30}\]
where  $\epsilon_i$  are complicated  homogeneous  polynomials in  the
$t_j's$  of degree  $i$ and  invariant  under the  permutation of  the
$t_j's.$ We define entire functions $\{P'_i(t)\}$ as follows,
\begin{equation}\label{A_n}
 A_n: P'_i=t_i-t_{i+1} \quad i=1,\cdots, n
\end{equation}
\begin{equation}\label{D_n}
D_n: P'_i=t_i-t_{i+1} \quad i=1,\cdots, n-1\quad \mbox{and}\quad P'_n=
t_{n-1}+t_{n}
\end{equation}
\begin{equation}\label{E_n} 
E_n:   P'_i=t_i-t_{i+1}\quad   i=1,\cdots,  n-1\quad   \mbox{and}\quad
P'_n=-t_1- t_2-t_3
\end{equation}

\section{ A description of N=1 ADE quivers and  reflection functors}
\label{section2}
\subsection{Describing N=1 ADE quivers}

This requires some detailed explanation, mainly of the relations
(\ref{relation1}), which distinguish these from ADE quivers.  To make
our presentation intelligible to non-experts, we briefly recall some
definitions and established facts. ( Here all vectors are over a field
$k.$)

A {\it quiver} $\Gamma=\left(V_\Gamma, E_\Gamma\right)$---without
relations---is a directed graph.

A {\it representation} $(V,f)$ of a quiver $\Gamma$ is an assignment
to each vertex $i\in V_\Gamma$ of a vector space $V(i),$ and to each
directed edge $ij\in E_\Gamma$ of a linear transformation $f_{ji}:
V(i)\to V(j).$

A $\mathit {morphism}$ $h:(V,f)\to (V',f')$ between 
  representations of $\Gamma$ over $k$ is a collection $\{h_i:V(i)\to
  V'(i)\}_{i\in V_\Gamma}$ of $k$-linear maps such that for each
  edge $ij\in E_\Gamma$ the diagram
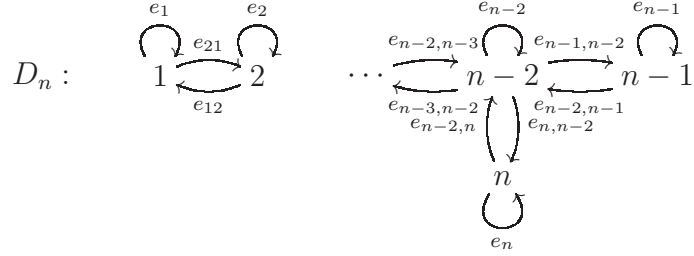
\[\begin{CD} V(i)@>h_i>> V'(i)\\
             @Vf_{ji} VV @Vf'_{ji} VV\\
             V(j) @>h_j>> V'(j)
\end{CD}\] 
commutes.  Compositions of morphisms are defined in the usual way. For
a path $p: i_1\to i_2\to \cdots \to i_r$ in $\Gamma,$ and a
representation $(V,f),$ we let $f_p$ be the composition of the linear
transformations $f_{i_{k+1}i_{k}}:V(i_k)\to V(i_{k+1}),$ $1\le k< r.$
And given vertices $i,\, j$ in $V_\Gamma,$ and paths $p_1,\cdots, p_n$
from $i$ to $j,$ a {\it relation} $\sigma$ on quiver $\Gamma$ is a
linear combination $\sigma=a_1p_1+\cdots +a_np_n,\, a_i\in k.$ If
$(V,f)$ is a representation of $\Gamma,$ we extend the $f$-notation by
setting $f_\sigma=a_1f_{p_1}+\cdots +a_nf_{p_n}: V(i)\to V(j).$ A {\it
  quiver with relations} is a pair $(\Gamma, \rho),$ where
$\rho=\left(\sigma_t\right)_{t\in T}$ is a set of relations on
$\Gamma.$ And a representation $(V,f)$ of $(\Gamma, \rho)$ is a
representation $(V,f)$ of $\Gamma$ for which $f_\sigma=0$ for all
relations $\sigma\in \rho.$ We can then define, in the obvious way,
{\it subrepresentations} $(V',f')$ of $(V,f),$ the {\it sum} of two
representations, and when a representation $(V,f)$ of $\left(\Gamma,
  \rho\right)$ is indecomposable, of finite representation type, and
simple.

\begin{definition}\label{N=1-A-D-E-Quiver}
  Given an $ADE$ Dynkin diagram $\Gamma,$
\[
\xymatrix{A_n: &
1 \edge{r} & 2 \edge{r} &  3\edge{r} &\cdots \edge{r}& n}
\]
\[
\xymatrix{& & & & n\\
D_n:& 1 \edge{r} & 2\edge{r} &\cdots & (n-2)\edge{l}\ar@{-}[u]\edge{r}& (n-1)}
\]
\[
\xymatrix{& & &  n\\
E_n:& 1 \ar@{-}[r] & 2\ar@{-}[r]& 3\ar@{-}[u]\ar@{-}[r] &\cdots & (n-1)\ar@{-}[l]}
\]
we can construct an associated directed graph, still denoted by
$\Gamma,$ by replacing each edge by a pair of oppositely directed
edges and associating to each vertex a loop.

\[
\xymatrix{  A_n:   &  1  \ar@(ul,ur)^{e_1}   \ar@/^/[r]^{e_{21}}  &  2
\ar@(ul,ur)^{e_2}   \ar@/^/  [l]^{e_{12}}   \ar@/^/[r]^{e_{32}}   &  3
\ar@(ul,ur)^{e_3}\ar@/^/[l]^{e_{23}} &  \cdots \ar@/^/[r] & \ar@/^/[l]
n \ar@(ul,ur)^{e_n}\\ }
\]
\[
\xymatrix{  D_n:   &  1  \ar@(ul,ur)^{e_1}   \ar@/^/[r]^{e_{21}}  &  2
\ar@(ul,ur)^{e_2}       \ar@/^/       [l]^{e_{12}}      &       \cdots
\ar@/^/[r]^{e_{n-2,n-3}}      &      n-2\ar@/^/[l]^      {e_{n-3,n-2}}
\ar@(ul,ur)^{e_{n-2}}\ar@/^/[d]^{e_{n,n-2}} \ar@/^/[r]^{e_{n-1, n-2}}&
n-1   \ar@(ul,ur)^{e_{n-1}}\ar@/^/[l]^{e_{n-2,n-1}}  \\   &  &   &  &n
\ar@(dl,dr)_{e_n} \ar@/^/[u]^{e_{n-2,n}} & }
\]
The associated directed graphs for $E_n$ ($n=6,\, 7,\, 8$) are,

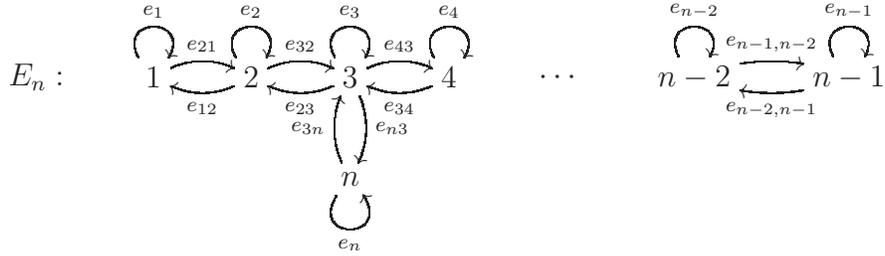
\[
\xymatrix{  E_n:   &  1  \ar@(ul,ur)^{e_1}   \ar@/^/[r]^{e_{21}}  &  2
\ar@(ul,ur)^{e_2}    \ar@/^/   [l]^{e_{12}}    \ar@/^/[r]^{e_{32}}   &
3\ar@(ul,ur)^{e_3}                                  \ar@/^/[r]^{e_{43}}
\ar@/^/[d]^{e_{n3}}\ar@/^/[l]^{e_{23}}      &     4\ar@/^/[l]^{e_{34}}
\ar@(ul,ur)^          {e_{4}}&          \cdots          &          n-2
\ar@(ul,ur)^{e_{n-2}}\ar@/^/[r]^{e_{n-1,n-2}}&       n-1      \ar@(ul,
ur)^{e_{n-1}}\ar@/^/[l]^{e_{n-2,n-1}}  \\  &  &  &n  \ar@(dl,dr)_{e_n}
\ar@/^/[u]^{e_{3n}} & }
\]
Then the $N=1$ ADE quivers are just the above associated directed
graphs, but with relations (\ref{relation1})

\begin{equation}\label{relation1}
\sum_is_{ij}e_{ji}e_{ij}+P'_j(e_j)=0, \quad
e_{ij}e_j=e_ie_{ij}
\end{equation}
where 
\begin{equation*}
\begin{cases}
s_{ij}=0 \quad \mbox{if}\quad i\quad \mbox{and}\quad j\quad \mbox{are not
adjacent}\\
s_{ij}= 1 \quad \mbox{if}\quad i\quad\mbox{and} \quad j\quad \mbox{are
adjacent and} \quad i>j\\
s_{ij}=- 1 \quad \mbox{if}\quad i\quad\mbox{and} \quad j\quad \mbox{are
adjacent and} \quad i<j\\
\end{cases}
\end{equation*}
and where $P'_j(x)$  is a certain fixed polynomial $\forall j.$  

If $(V,f)$ is a representation of an $N=1$ ADE quiver, the
corresponding structures are
 
\[
\xymatrix{  A_n:  &V(1)  \ar@(ul,ur)^{\Phi_1} \ar@<1ex>[r]^{Q_{21}}  &
V(2) \ar@(ul,ur)^ {\Phi_2} \ar@<1ex>[l]^{Q_{12}} \ar@<1ex>[r]^{Q_{32}}
&    V(3)   \ar@(ul,ur)^   {\Phi_3}\ar@<1ex>[l]^{Q_{23}}    &   \cdots
\ar@<1ex>[r] & \ar@<1ex>[l] V(n) \ar@(ul,ur)^{\Phi_n}\\ }
\]
\[
\xymatrix{  D_n:  &V(1)  \ar@(ul,ur)^{\Phi_1} \ar@<1ex>[r]^{Q_{21}}  &
V(2)    \ar@(ul,ur)^   {\Phi_2}    \ar@<1ex>[l]^{Q_{12}}    &   \cdots
\ar@<1ex>[r]^{Q_{n-2,n-3}}   &  V(n-2)  \ar@<1ex>[l]^(.8){Q_{n-3,n-2}}
\ar@(ul,ur)^{\Phi_{n-2}}\ar@<1ex>[d]^(.5)                   {Q_{n,n-2}}
\ar@<1ex>[r]^{Q_{n-1,n-2}}&       V(n-1)      \ar@(ul,ur)^{\Phi_{n-1}}
\ar@<1ex>[l]^{Q_{n-2,n-1}}  \\   &  &  &   &V(n)  \ar@(dl,dr)_{\Phi_n}
\ar@<1ex>[u]^{Q_{n-2,n}} & }
\]
\[
\xymatrix{  E_n: &  V(1) \ar@(ul,ur)^{\Phi_1}  \ar@<1ex>[r]^{Q_{21}} &
V(2) \ar@(ul,ur)^ {\Phi_2} \ar@<1ex>[l]^{Q_{12}} \ar@<1ex>[r]^{Q_{32}}
&                       V(3)\ar@(ul,ur)^                      {\Phi_3}
\ar@<1ex>[r]^{Q_{43}}\ar@<1ex>[d]^{Q_{n3}}\ar@<1ex>[l]^{Q_{23}}       &
V(4)\ar@<1ex>[l]^{Q_{34}}  \ar@(ul,ur)^{\Phi_{4}}&   \cdots  &  V(n-2)
\ar@(ul,ur)^{\Phi_{n-2}}\ar@<1ex>[r]^{Q_{n-1,n-2}}                    &
V(n-1)\ar@(ul,ur)^ {\Phi_{n-1}}\ar@<1ex>[l]^{Q_{n-2,n-1}} \\ & & &V(n)
\ar@(dl,dr)_{\Phi_n} \ar@<1ex>[u]^{Q_{3n}} & }
\]
where we have write $Q_{ij}=f_{e_{ij}}, \, \Phi_j=f_{e_j}.$  And the relation (\ref{relation1}) becomes
\begin{equation}\label{relation}
\sum_is_{ij}Q_{ji}Q_{ij}+P'_j(\Phi_j)=0, \quad
Q_{ij}\Phi_j=\Phi_iQ_{ij}.
\end{equation}
\end{definition}

\subsection{Reflection functors}\label{defreflec}
Given an $N=1$ ADE quiver $\Gamma$ and $k\in V_\Gamma,$  denote by
$\Gamma^+_k$ the quiver defined by dropping all arrows starting from $k,$
and denote by $\Gamma^-_k$ the quiver defined from $\Gamma$ by dropping
all arrows ending at $k.$

Given a representation $V$ of an $N=1$ ADE quiver $\Gamma,$  we can define
a representation of $\Gamma^+_k$ which we still denote  as $V,$   by
forgetting all maps which have domain $V(k).$  Similarly, we define a
representation of $\Gamma^-_k$ which we still denote  by $V,$   by
forgetting all maps which has range $V(k).$  Then we can apply the
reflection functor $F^+_k$ in \cite{BGP} to the representation $V$ of
$\Gamma^+_k$ and   apply the reflection functor $F^-_k$ in \cite{BGP} to
the representation $V$ of $\Gamma^-_k.$

In the following definition \ref{reffunctor}, we modify the reflection
functors in \cite{BGP} for the purpose of this paper.
\begin{definition}\label{reffunctor}Let $\Gamma$ be an $N=1$ ADE quiver and 
$k$ a vertex of $\Gamma.$  Let
  $$\Gamma^k=\{i\, \mid \, i \mbox{ adjacent to}\, k\}$$
  
For a quiver representation $W$ of $\Gamma^+_k,$  define a representation 
$F^{+}_k(W)$ of $\Gamma^{-}_k$ by 
\begin{equation}
F^{+}_k(W)(i)=\begin{cases} W(i) \quad \mbox{if}\quad i\not= k\\
                               \ker h \quad \mbox{if}\quad i=k\\
\end{cases}
\end{equation}
where
\[h:\oplus_{i\in \Gamma^k} W(i)\to W(k)\]
is defined by
\[h\left((x_i)_{i\in \Gamma^k}\right)=\sum_{i\in \Gamma^k}s_{ik}Q_{ki}x_i\]
If  $i, \, j\not=k,$ we  define $Q'_{ij}=Q_{ij}: W(j)\to W(i).$  
 If $i\in\Gamma^k,$ define $Q'_{ik}: F^{+}_k(W)(k)\to W(i)$ by 
\begin{equation}\label{twistequation}
Q'_{ik}(x_j)_{j\in \Gamma^k}=-s_{ki}x_i
\end{equation} 
For a quiver representation $U$ of $\Gamma^-_k,$ define a representation 
$F^{-}_k(U)$ of $\Gamma^+_k$ by
\begin{equation}
F^{-}_k(U)(i)=\begin{cases} U(i) \quad \mbox{if}\quad i\not= k\\
                               \operatorname{coker}g \quad \mbox{if}\quad 
i=k\\
\end{cases}
\end{equation}
where 
\[g:U(k)\to \oplus_{i\in \Gamma^k} U(i)\]
is defined by
\[g(x)=(Q_{ik}x)_{i\in \Gamma^k}\]
and define $Q'_{ki}: U(i)\to F^{-}_k(U)(k)$ by  the natural  composition 
of
\begin{equation}
U(i)\to \oplus_{j\in \Gamma^k}U(j)\to F^{-}_k(U)(k)
\end{equation}
\end{definition}

\begin{remark} The definitions of the $F^{+}_k(W)$ and $Q'_{ik}$  in
\ref{reffunctor} are different than the corresponding definitions in
\cite{BGP}, while  $F^{-}_k(U)$ and $Q'_{ki}$  in \ref{reffunctor} are the
same as the corresponding definitions in \cite{BGP}.
\end{remark}

\subsection{The action of the  Weyl group on $\{P'_i\}, \, 1\le i \le n$}
 By \cite{KM92}, pp 461 and
463, we know the Weyl group $\mathfrak{W}_{A_n}$ of $A_n$ is generated by
reflections $r_1,\cdots, r_n$, which act as permutations of $t_1, \cdots,
t_{n+1},$  where $t_i$ is defined on Section~\ref{fibration}.

 In the $A_n$ case,  we can write $P'_i(x)$ in the relation  given in
(\ref{relation}) as
\[{A_{n}: \quad P'_i = t_i - t_{i+1}  \quad  i=1, \ldots ,n}\]
Then the action of $r_i$ extends to polynomials in the $\{t_i\}.$
By \cite{KM92}, pp 461 and 463, we know the Weyl group
$\mathfrak{W}_{D_n}$ of $D_n$  is generated by reflections $r_i$, for
$1\le i\le n-1,$  which act as permutations of $t_1,\cdots, t_n,$
together with $r_n$ which is defined by $$r_n(t_i)=\begin{cases}t_1 &
\text{if $1\le i\le n-2$}\\
                            -t_n & \text{if  $i=n-1$}\\
                            -t_{n-1} & \text {if $i=n$}\\ 
               \end{cases}$$
 In the $D_n$ case, we can write $P'_i(x)$ in the relation  given in
(\ref{relation}) as
$$D_{n}: \quad P'_{i}= t_{i} - t_{i+1}  \quad  i=1,\ldots ,n-1$$
and   $$P'_n=t_{n-1}+t_n$$
Again, the action of $r_i$ extends to polynomials in the $\{t_i\}.$
{By \cite{KM92}, pp 461 and 463, we know that the Weyl group
$\mathfrak{W}_{E_n}$ of $E_n$ is generated by reflections $r_i$ for $1\le
i\le n-1,$ which act as permutations of $t_1,\cdots , t_n,$  together with
$r_n,$ which is defined by
$$r_n(t_i)=\begin{cases}t_i-\frac{2}{3}(t_1+t_2+t_3) \quad \text{if} \quad
1\le i\le 3 \\ t_i+\frac{1}{3}(t_1+t_2+t_3) \quad \text{if} \quad 4\le
i\le n \\ \end{cases}$$  In the $E_n$ case, we can write $P'_i(x)$ in the
relation  given in (\ref{relation}) as
\[{E_{n}: \quad P'_i = t_i - t_{i+1}  \quad  i=1, \ldots ,n} \quad
\text{and} \quad P'_{n}=-t_1-t_2-t_3\]
Once again, the action of $r_i$ extends to polynomials in the $\{t_i\}.$
Based on these definitions of $r_i,\, 1\le i\le n,$  one can easily get
the following Lemma~\ref{weylgroup}.
\begin{lemma}\label{weylgroup} Let $\mathfrak{W}_\Gamma$ be the Weyl group
of the Dynkin diagram $\Gamma$ and let  $r_i\in \mathfrak{W}_\Gamma$
($1\le i\le n$) be a set of generators of reflections.   If $j$ is
distinct from $i$ and  not adjacent to $i,$ then
$r_i(P'_j(\Phi_j))=P'_j(\Phi_j).   $  If $j$ is adjacent to $i$ and
$j\not=i,$  then $r_i(P'_j(\Phi_j))=P'_j(\Phi_j)+P'_i(\Phi_j).$  Finally,
$r_i(P'_i(\Phi_i))=-P'_i(\Phi_i).$
\end{lemma}
\section{finite-to-one correspondence}\label{section3}
 In this section, we will give a   proof, using reflection functors, that
in the case of simple and distinct roots, the irreducible quiver
representations are in finite-to-one correspondence with the contractible
curves in the threefold.

\subsection{Applying the reflection functors to $N=1$ ADE quiver
representations}
Let $\Gamma$ be an N=1  ADE quiver.  Let
\[
{\mathcal A}_\Gamma=\left\{\sum_i n_i P_i' | n_i \in {\mathbb Z},\, \mbox{not
all}\quad n_i\quad \mbox{zero}\right\}
\]
where $P'_i,\, 1\le i\le n,$ are the polynomials in relation
(\ref{relation})

(*) Suppose no two elements $\sum n_i P_i',\, \sum m_i P_i'$ of the set
${\mathcal A}_\Gamma$
have a common root unless there is a constant c with $m_i = c n_i$ for all
$i.$
\begin{lemma}$(*)$ holds for any generic collection of polynomials
$P_i'$ of positive degree.
\end{lemma}
\begin{proof}
Left to the reader.
\end{proof}

\begin{lemma}\label{remark1} Let $V$ be an $N=1$ ADE quiver
representation,let $v_j$ be a $\lambda$- eigenvector of $\Phi_j$.  Then
$Q_{ij}\Phi_jv_j$ is either a $\lambda$-eigenvector of $\Phi_i$ or 0.
\end{lemma}
\begin{proof}
If $v_j$ is an eigenvector of $\Phi_j$ corresponding to eigenvalue
$\lambda,$ then from (\ref{relation}), we get
\[Q_{ij}\Phi_jv_j=\Phi_iQ_{ij}v_j\]
which implies that
\[\lambda Q_{ij}v_j=\Phi_iQ_{ij}v_j\]
Hence, $Q_{ij}v_j$ is either an eigenvector of $\Phi_i$ corresponding to
eigenvalue $\lambda$ or a $0$ vector.
\end{proof}
\begin{lemma}Let $V$ be a simple  representation of an $N=1$ ADE quiver.
Then there exists $\lambda$ such that if $v_i\in V(i)\not= 0,$ then
$\Phi_i v_i=\lambda v_i.$
\end{lemma}
\begin{proof} Let $\mathcal{A}=\{d| V(d)\not=0\}.$  Then $\mathcal{A}$ is
connected.  Otherwise, $V$ is not simple.    Let $a=\min \mathcal{A},$
then $\Phi_a$ has a eigenvector $v_a$ with eigenvalue $\lambda$.    For
$l\in \mathcal{A},$  let  $U(l)$ be the $\lambda$-eigenvector space of
$\Phi_l$.   By Lemma~\ref{remark1},  it's easy to see that $(W,g)=\{U(l):
l\in \mathcal{A}\}$ is a sub-representation of $V.$   Since $V$ is simple,
$(W,g)=V,$  which proves the result.
\end{proof}
Therefore, to show that we have only finitely many  simple
representations, it suffices to consider  representations $V$ for which
there exists a $\lambda$ such that if  $0\not= v_d\in V(d),$ then  $\Phi_d
v_d=\lambda v_d.$
{\it  In the rest of this section, we only consider  quiver
representations $V$ with this property.}

\begin{lemma}\label{F^+} Let $V$ be a simple  representation of an $N=1$
  ADE quiver.  Suppose $V$ is not concentrated at vertex $k$.  Then
\[
\dim\left(F^{+}_k(V)\right)_k=\sum_{i\in\Gamma^k}\dim\left(V(i)\right)-\dim\left(V(k)\right)
\]

\end{lemma}
\begin{proof} We know that $\left( F^{+}_k(V)\right)(k)=\ker h,$  where
$h:\oplus_{i\in \Gamma^k} V(i)\to V(k)$  is defined by
\[h(x_i)_{i\in \Gamma^k}=\sum_{i\in \Gamma^k}s_{ik}Q_{ki}x_i\]

To prove the lemma is equivalent to prove that $h$ is surjective.

Case I: $V(k)\not=0.$ If $h$ is not surjective and $h\not=0,$ then we
can replace $V(k)$ by $h\left(\oplus_{i\in \Gamma^k}V(i)\right)$ to
get a sub-representation of $V.$ But this contradicts  the
simplicity of $V.$ Case II: $V(k)=0.$ We get that $h$ is surjective
since $h\equiv 0$ in this case.
\end{proof}
\begin{lemma}\label{F^-} Let $V$ be a simple  representation of an $N=1$
ADE quiver.   Suppose $V$ is not  concentrated at vertex $k$. Then
\[
\dim\left(F^{-}_k(V)\right)_k=\sum_{i\in\Gamma^k}\dim\left(V(i)\right)-\dim\left(V(k)\right)
\]
\end{lemma}
\begin{proof} We know that $\left(F^{-}_k(V)\right)(k)=\coker g,$  where
$g:V(k)\to \oplus_{i\in \Gamma^k}V(i)$  is defined by
$g(x)=\left(Q_{ik}x\right)_{i\in \Gamma^k}.$ 
To prove the lemma is equivalent to prove that $g$ is injective.

Case I: $V(k)\not=0.$
 If $\ker g\not=0,$  then we
can define a simple sub-representation which concentrated at vertex $k.$
This contradicts  the simplicity of $V.$ 

Case II: $V(k)=0.$
We get that $g$ is injective since $g\equiv 0$ in this case.
 \end{proof}

\begin{lemma}\label{iso} Let $V$ be a simple  representation of an $N=1$
  ADE quiver.  Suppose $V$ is not concentrated at vertex $k$. If
  $P'_k(\lambda)\ne 0,$ then there is a natural isomorphism $\varphi$
  between $F^{+}_k(V)(k)$ and $F^{-}_k(V)(k).$
\end{lemma}
\begin{proof}
Let $g: V(k)\to \oplus_{i\in \Gamma^k}V(i)$ be defined by
\[
g(x)=\left(Q_{ik}x\right)_{i\in \Gamma^k}
\]
and $h: \oplus_{i\in \Gamma^k} V(i)\to V(k)$ be defined by

\[
h(x_i)_{i\in \Gamma^k}=\sum_{i\in \Gamma^k}s_{ik}Q_{ki}x_i
\]

We have
\[F^{-}_k(V)(k)=\coker g\]
and
\[F^{+}_k(V)(k)=\ker h\]
Since $V$ is simple and not concentrated at $k,$ $g$ is injective and
$h$ is surjective.  We have
\begin{eqnarray*}
 \dim F^{+}_k(V)(k)
=\dim F^{-}_k(V)(k)
=\sum_{i\in \Gamma^k} \dim V(i)-\dim V(k)
\end{eqnarray*}
\[
\xymatrix{ &&  0\ar[d] \\ && F^{+}_k(V)(k)\ar[d]^{h'}\ar[rd]^{\varphi}\\ 0\ar[r] &
V(k)\ar[rd]_{-P'_k(\lambda)}\ar[r]^{g}& \oplus_{i\in
\Gamma^k}V(i)\ar[d]^{h}\ar[r]^{g'}& F^{-}_k(V)(k)\ar[r]& 0\\
&& V(k)\ar[d]\\ && 0\\}
\]
Since $P'_k(\lambda)\ne 0,$ $\im g\cap F^{+}_k(V)(k)=\{0\}.$ Let $g':
\oplus_{i\in \Gamma^k}V(i)\to F^{-}_k(V)(k)$ be the natural surjective
map induced by $g$ and let $h': F^{+}_k(V)(k)\to \oplus_{i\in
  \Gamma^k}V(i)$ be the natural inclusion map induced by $h.$ Then
$\varphi=g'\circ h': F^{+}_k(V)(k)\to F^{-}_k(V)(k)$ is a natural
isomorphism (Since $\dim F^{+}_k(V)(k)=\dim F^{-}_k(V)(k)$ and
$\varphi$ is injective by $\im g\cap F^{+}_k(V)(k)=\{0\}.$)
\end{proof}
\begin{definition}\label{newrep} By Lemma~\ref{iso}, if $P'_k(\lambda)\ne 0,$  we can
  construct a new representation $F_k(V)$ of $\Gamma$ by
\[
F_k(V)(i)=\begin{cases} V(i) \quad \mbox{if} \quad i\not= k\\
             F^+_k(V)(k)  \quad \mbox{if} \quad i=k\\
\end{cases}
\]
defining $Q'_{lk}$ as it is defined in $F^{+}_k(V)(k)$ and defining
$Q'_{km}$ as the composition map $P'_k(\lambda)\cdot\varphi^{-1}\circ
{\underline {Q'_{km}}}: V(m)\to F^{+}_k(V)(k),$ where $\underline
{Q'_{km}}:V(m)\to F^{-}_k(V)(k)$ is the natural map defined in
$F^{-}_k(V)$ and $\varphi: F^{+}_k(V)(k)\to F^{-}_k(V)(k)$ is the
isomorphism defined in Lemma~\ref{iso}.

If $V$ is simple, we define 
\[
\Phi'_i: F_k(V)(i)\to F_k(V)(i)
\]
by $\Phi'_i(x)=\lambda x,$ where $\lambda$ is the eigenvalue of $\Phi$
on $V(i)$ that appeared in the representation of $V.$ Abusing
notation, we still denote $\Phi'_i$ as $\Phi_i.$
\end{definition}

\begin{lemma}\label{poly1}If $V$ is a simple  representation   of  $N=1$
ADE quiver, then
\[\sum_i\dim(V(i))\cdot P'_i(\lambda)=0\]
\end{lemma}
\begin{proof}This follows from the fact that $\forall$ pair $i$ and $j,$
$\Tr \,Q_{ij}Q_{ji}=\Tr\, Q_{ji}Q_{ij},$  and $\forall, \, i,$ $\Tr\,
\Phi_i=\lambda\cdot \dim V(i),$  where $\lambda$ is an eigenvalue for all
$\Phi_i.$
Now take trace operation to relations (\ref{relation}) and then sum the
resulting equations.   The result follows.
\end{proof}
\begin{lemma}\label{poly2} Let $V$ be a simple representation of an $N=1$
ADE quiver $\Gamma,$  not concentrated at vertex $k$.  Then
\[\sum\dim \left(F_k(V)\right)(i) r_k\left(P'_i(\lambda)\right)=\sum\dim
V(i) P'_i(\lambda)\]
\end{lemma}
\begin{proof}\begin{eqnarray*}
&&\sum\dim \left(F_k(V)\right)(i) r_k\left(P'_i(\lambda)\right)\\
&=&\sum_{i\in \Gamma^k}\dim V(i)\left(P'_k(\lambda)+P'_i(\lambda)\right)
+\sum_{j\in \Gamma-\Gamma^k}\dim V(j)P'_j(\lambda)\\
&+&\left(-\dim V(k)+\sum_{i\in \Gamma^k}\dim
V(i)\right)\left(-P'_k(\lambda)\right)\\
&=&\sum \dim V(i)P'_i(\lambda)\\
\end{eqnarray*}
\end{proof}

\begin{proposition}\label{newrelation} Let $V$ be a simple representation
  of an $N=1$ ADE quiver which is not concentrated at vertex $k.$ If
  $P'_k(\lambda)\ne 0,$ then $F_k(V)$ satisfies the following new
  relations
\begin{equation}
\sum_is_{ij}Q'_{ji}Q'_{ij}+r_k\left(P'_j(\Phi_j)\right)=0, \quad
Q'_{ij}\Phi_j=\Phi_iQ'_{ij}.
\end{equation}
\end{proposition}
\begin{proof}

If $i\not\in \Gamma^k$ and $i\not=k,$ where $i$ is a vertex of
$\Gamma$ such that $V(i)\ne 0,$ there is nothing to prove. For $ j\in
\Gamma^k\cup \{k\},\, b\in \left(F_k(V)\right)_j,$ we have
\[Q'_{ij}\Phi_j b=\lambda Q'_{ij}b=\Phi_iQ'_{ij}b\]

For $i\in \Gamma^k$ and $x\in V(i)$, by Definition~\ref{newrep}, we know that
\[
Q'_{ki}x=P'_k(\lambda)\cdot {\varphi}^{-1}\circ\underline{Q'_{ki}}x
\]
where
 $\underline{Q'_{ki}}x=[(x_j)_{j\in
\Gamma^k}]\in F^{-}_k(V)(k),$ and
\[
x_j=\begin{cases} 0 \quad \mbox{if} \quad j\not=i\\
                  x  \quad \mbox{if} \quad j=i\\
\end{cases}
\]
After a short computation, we see that
\[
Q'_{ki}x=(y_j)_{j\in \Gamma^k}
\]

where
\[
y_j=\begin{cases} P'_k(\lambda)x+s_{ik}Q_{ik}Q_{ki}x\quad  \mbox{if} \quad j=i\\
Q_{jk}s_{ik}Q_{ki}x \quad \mbox{if} \quad j\not= i\\
\end{cases}
\]

It follows that
\begin{equation*}
s_{ki}Q'_{ik}Q'_{ki}x = s_{ki}Q'_{ik}(y_j)_{j\in
\Gamma^k}=-P'_k(\lambda)x-Q_{ik}s_{ik}Q_{ki}x
\end{equation*}

Hence for $i\in \Gamma^k$  we have
\begin{eqnarray*}
&& \sum_{j}s_{ji}Q'_{ij}Q'_{ji}x+r_k\left(P'_i(\lambda)\right)x\\
&=&\sum_{j}s_{ji}Q'_{ij}Q'_{ji}x+r_k\left(P'_i(\lambda)\right)x\\
&=&\sum_{j}s_{ji}Q'_{ij}Q'_{ji}x+P'_i(\lambda)x+P'_k(\lambda)x\\
&=&\sum_{ j\not=k}s_{ji}Q_{ij}Q_{ji}x
+s_{ki}Q'_{ik}Q'_{ki}x+P'_i(\lambda)x+P'_k(\lambda)x\\
&=&\sum_{
j\not=k}s_{ji}Q_{ij}Q_{ji}x-P'_k(\lambda)x-Q_{ik}s_{ik}Q_{ki}x+P'_i(\lambda)x+P'_k(\lambda)x\\
&=& 0\\
\end{eqnarray*}
Let $(x_i)_{i\in \Gamma^k}\in F^{+}_k(V)(k).$  Then
\[
s_{ik}Q'_{ki}Q'_{ik}(x_i)_{i\in \Gamma^k}=Q'_{ki}x_i=\left(x_{i_j}\right)_{j\in
\Gamma^k}
\]

where
\[
x_{i_j}=\begin{cases} P'_k(\lambda)x_i+Q_{ik}s_{ik}Q_{ki}x_i \quad \mbox{if}
\quad j=i\\
Q_{jk}s_{ik}Q_{ki}x_i \quad \mbox{if} \quad j\not=i\\
\end{cases}
\]
Hence we have
\begin{eqnarray*}
&&\sum_{i\in \Gamma^k}s_{ik}Q'_{ki}Q'_{ik}(x_i)_{i\in
\Gamma^k}+r_k(P'_k(\lambda))(x_i)_{i\in \Gamma^k}\\
&=&\sum_{i\in \Gamma^k}s_{ik}Q'_{ki}Q'_{ik}(x_i)_{i\in
\Gamma^k}-P'_k(\lambda)(x_i)_{i\in \Gamma^k}\\
&=& \sum_{i\in \Gamma^k}\left(x_{i_j}\right)_{j\in
\Gamma^k}-P'_k(\lambda)(x_i)_{i\in \Gamma^k}\\
&=&0\\
\end{eqnarray*}

\end{proof}

\begin{lemma}\label{ref1} If $V$ is a simple representation of  an $N=1$ ADE 
  quiver which is not concentrated at vertex $k$ and if $(*)$ holds,
  then $F_kF_k(V)\cong V.$ Consequently, $F_k(V)$ is a simple
  representation.
\end{lemma}
\begin{proof}
  We know that $Q'_{ki}: V(i)\to F_k(V)(k)$ is defined by
  $Q'_{ki}x_i=P'_k(\lambda)\varphi^{-1}\underline{Q_{ki}}x_i,$ where
  $\underline{Q_{ki}}: V(i)\to F^{-}_k(V)(k)$ is the composition of
  $V(i)\to \oplus_{i\in \Gamma^k}V(i)$ and $\oplus_{i\in
    \Gamma^k}V(i)\to F^{-}_k(V)(k)$ (See Definition~\ref{newrep}).  We
  also know that
\[
F_kF_k(V)(k)=\left\{(x_i)\in \oplus_{i\in \Gamma^k} V(i) \,\mid\, \sum_{i\in \Gamma^k}s_{ik}Q'_{ki}x_i=0\right\}
\]
We have
\[
\sum_{i\in \Gamma^k}s_{ik}Q'_{ki}x_i=P'_k(\lambda)\varphi^{-1}\sum_{i\in \Gamma^k}s_{ik}\underline{Q_{ki}}x_i
\]
Since $P'_k(\lambda)\not=0$ and $\varphi$ is an isomorphism, we get
\[
F_kF_k(V)(k)=\left\{(-s_{ki}Q_{ik}x)\,\mid\, x\in V(k)\right\}
\]

Let $g:V\to F_kF_k(V)$ be defined in the following way:
\[
g_i=\begin{cases} i: V(i)\to F_kF_k(V)(i)=V(i) \quad \mbox{if}\, i\not=k\\
(-s_{ki}Q_{ik})\quad \mbox{if} \, i=k\\
\end{cases}
\]
where $i:V(i)\to F_kF_k(V)(i)=V(i)$ is the identity map.

Then it is clear  that (\ref{[1]}) is commutative.

\begin{equation}\label{[1]}
\xymatrix{ V(k)\ar[r]^{Q_{ik}}\ar[d]^{g_k} & V(i)\ar[d]^{g_i}\\ F_kF_k(V)(k)\ar[r]^{Q''_{ik}}& V(i)\\}
\end{equation}
Let's check the commutativity of (\ref{[2]}).

\begin{equation}\label{[2]}
\xymatrix{V(i)\ar[r]^{Q_{ki}}\ar[d]^{g_i}& V(k)\ar[d]^{g_k}\\ V(i)\ar[r]^{Q''_{ki}}& F_kF_k(V)(k)\\}
\end{equation}
Let $(Q''_{ki}x_i)_j$ ( resp. $(Q'_{ki}x_i)_j$)  denote the $j-$th coordinate of $Q''_{ki}x_i$ (resp. $Q'_{ki}x_i$).
We know that
\[
(Q''_{ki}x_i)_j=\begin{cases} -P'_k(\lambda)x_i+Q'_{ik}s_{ik}Q'_{ki}x_i\\
                          Q'_{jk}s_{ik}Q'_{ki}x_i\quad \mbox{if} \, j\not=i\\
              \end{cases}
\]
where  
\[
(Q'_{ki}x_i)_j=\begin{cases} P'_k(\lambda)x_i+Q_{ik}s_{ik}Q_{ki}x_i\\
                          Q_{jk}s_{ik}Q_{ki}x_i\quad \mbox{if} \, j\not=i\\
              \end{cases}
\]
If  $i>k,$  then we have
\begin{eqnarray*}
(Q''_{ki}x_i)_i&=&-P'_k(\lambda)x_i+Q'_{ik}s_{ik}Q'_{ki}x_i\\
&=&-P'_k(\lambda)x_i+P'_k(\lambda)x_i+Q_{ik}s_{ik}Q_{ki}x_i\\
&=&Q_{ik}s_{ik}Q_{ki}x_i=Q_{ik}Q_{ki}x_i\\
&=&Q_{ik}Q_{ki}x_i= -s_{ki}Q_{ik}Q_{ki}x_i\\
\end{eqnarray*}
If $i>k$ and  $j>k$, then we have
\[
(Q''_{ki}x_i)_j=Q'_{jk}s_{ik}Q'_{ki}x_i=Q_{jk}s_{ik}Q_{ki}x_i=Q_{jk}Q_{ki}x_i=-s_{kj}Q_{jk}Q_{ki}x_i
\]
If $i>k$ and $j<k,$ then we have
\[
(Q''_{ki}x_i)_j=Q'_{jk}s_{ik}Q'_{ki}x_i=Q'_{jk}Q'_{ki}x_i=-Q_{jk}Q_{ki}x_i=-s_{kj}Q_{jk}Q_{ki}x_i
\]
If $i<k,$  then we have
\begin{eqnarray*}
(Q''_{ki}x_i)_i&=&-P'_k(\lambda)x_i+Q'_{ik}h'_{ki}x_i\\
&=&-P'_k(\lambda)x_i-Q'_{ik}Q'_{ki}x_i\\
&=&-P'_k(\lambda)x_i+(P'_k(\lambda)x_i+Q_{ik}h_{ki}x_i)\\
&=&Q_{ik}h_{ki}x_i=-Q_{ik}Q_{ki}x_i\\
&=&-s_{ki}Q_{ik}Q_{ki}x_i\\
\end{eqnarray*}
If $i<k$ and $j>k,$ then we have
\[
(Q''_{ki}x_i)_j=Q'_{jk}s_{ik}Q'_{ki}x_i=-Q'_{jk}Q'_{ki}x_i=-Q_{jk}s_{ik}Q_{ki}x_i=Q_{jk}Q_{ki}x_i=-s_{kj}Q_{jk}Q_{ki}x_i
\]
If $i<k$ and $j<k,$ then we have
\[
(Q''_{ki}x_i)_j=Q'_{jk}s_{ik}Q'_{ki}x_i=-Q'_{jk}Q'_{ki}x_i=Q_{jk}s_{ik}Q_{ki}x_i=-Q_{jk}Q_{ki}x_i=-s_{kj}Q_{jk}Q_{ki}x_i
\]
Therefore the diagram (\ref{[2]}) is commutative.

  Diagram (\ref{[3]}) is commutative
since $\lambda$ is a common eigenvalue of $V(k)$ and $F_kF_k(V)(k).$

\begin{equation}\label{[3]}
\xymatrix{V(k)\ar[d]^{g_k}\ar[r]^{\lambda}&      V(k)\ar[d]^{g_k}\\      F_kF_k(V)(k)\ar[r]^{\lambda}&
          F_kF_k(V)(k)\\}
\end{equation}

We prove the later part of the Lemma here. Since $V$ is not
concentrated at vertex $k,$ $\exists m\not= k$ such that $V(m)\not=0.$
It follows that $F_k(V)(m)=V(m)\not=0.$  Let $v\in F_k(V)(m)$ be an
nonzero element. If $F_k(V)$ is not simple, then there exists,
actually we can construct a simple sub-representation $W$ of $F_k(V),$
such that $v\in W(m).$  It follows that $F_k(W)$ is a proper
sub-representation of $F_kF_k(V).$ Since $F_kF_k(V)\cong V$ and $V$ is
simple, this leads a contradiction.
   
\end{proof}

\begin{corollary}\label{cor} Let $\Gamma$ be an $N=1$ ADE quiver.  Let
\[
\mathcal {B}_\Gamma=\{r_i  (P'_j(x))  | r_i  \in \mathfrak {W}_\Gamma\}
\]
where $\mathfrak {W}_\Gamma$ is the Weyl group of $\Gamma$ and $P'_j$
is the polynomial defined on relation~(\ref{relation}).  Suppose $(*)$
holds and each element in $\mathcal {B}_\Gamma$ has simple roots.  If
$V$ is a simple representation, then either $F_k(V)$ is simple or
$V\cong L_k,$ where $L_k$ is a simple representation concentrated at
vertex $k.$
\end{corollary}
\begin{proof}
  If $V$ is simple and concentrated at vertex $k,$ then $V\cong L_k,$
  where $L_k$ is a simple representation concentrated at vertex $k.$
  Assume $V$ is not concentrated at vertex $k.$ Since $V$ is simple,
  by Lemma~\ref{F^+} and Lemma~\ref{F^-}, we can apply $F_k$ to $V.$
  Then $F_k(V)$ is simple by the later part of Lemma~\ref{ref1}.
\end{proof}

\subsection {A proof of the Main Theorem}\label{applyref}

Let $\Gamma$ be a quiver. Following \cite{BGP}, for a representation
$V,$ we define $\dim (V)=\left(\dim V(i)\right)_{i\in V_\Gamma}.$
Denote by ${\mathcal C}_\Gamma=\{x=(x_\alpha)\mid x_\alpha\in \Q,\,
\alpha\in V_\Gamma\},$ where $\Q$ denotes the set of rational numbers.
We call a vector $x=(x_\alpha)$ positive (written $x>0$) if $x\not= 0$
and $x_\alpha\ge 0$ for all $\alpha\in V_\Gamma.$ For each $\beta\in
V_\Gamma,$ denote by $\sigma_\beta$ the linear transformation in
${\mathcal C}_\Gamma $ defined by the formula $(\sigma_\beta
x)_\gamma=x_\gamma$ for $\gamma\not= \beta,$ $(\sigma_\beta
x)_\beta=-x_\beta+\sum_{l\in \Gamma^\beta}x_l,$ where $l\in
\Gamma^\beta$ is the set of vertices adjacent to $\beta.$

For each vertex $\alpha\in V_\Gamma$ we denoted by $\Gamma_\alpha$ the
set of edges containing $\alpha.$ Let $\Lambda$ be an orientation of
the graph $\Gamma.$ We denote by $\sigma_{\alpha}\Lambda$ the
orientation obtained from $\Lambda$ by changing the directions of all
edges $l\in \Gamma_\alpha.$ Following \cite{BGP}, we say that a vertex
$i$ of a quiver $\Gamma$ with orientation $\Lambda$ is (-)-accessible
(resp.  (+)-accessible) if for any edge $e$ having $i$ as a vertex, we
have the final vertex of $f(e)$ of $e$ satisfying $f(e)\not=i$ (resp.
the initial vertex $i(e)$ of $e$ satisfying $i(e)\not=i.$) We say that
a sequence of vertices $\alpha_1,\, \alpha_2,\cdots, \alpha_k$ is
(+)-accessible with respect to $\Lambda$ if $\alpha_1$ is
(+)-accessible with respect to $\Lambda,$ $\alpha_2$ is (+)-accessible
with respect to $\sigma_{\alpha_1}\Lambda, $ $\alpha_3$ is
(+)-accessible with respect to
$\sigma_{\alpha_2}\sigma_{\alpha_1}\Lambda,$ and so on.  We define a
(-)accessible sequence similarly.

\begin{definition} Let $\Gamma$ be a graph without loops. We denote by  
$\mathscr{C}_\Gamma$ the linear space over $\mathbb{Q}$ consisting of sets 
$x=(x_\alpha)$ of rational numbers $x_\alpha$ ($\alpha\in \Gamma_V$).
  We call a vector $x=(x_\alpha)$ {\it positive} (written $x>0$) if
  $x\not=0$ and $x_\alpha\ge 0$ for all $\alpha\in \Gamma_V.$ We
  denote by $B$ the quadratic form on the space $\mathscr{C}_\Gamma$
  defined by the formula $B(x)=\sum x^2_\alpha -\sum_{l\in
    \E_\Gamma}x_{r_1(l)}x_{r_2(l)},$ where $r_1(l)$ and $r_2(l)$ are
  the ends of the edge $l.$ We denote by $<,>$ the corresponding
  symmetric bilinear form.
\end{definition}
\begin{lemma}\cite[Lemma 2.3]{BGP} Suppose that the form $B$ for the graph $\Gamma$ 
is positive definite. Let $c=\sigma_n\cdots \sigma_2\sigma_1.$  If $x\in \mathscr{C}_\Gamma,$  
$x\not=0,$  then for some $i$ the vector $c^ix$ is not positive.
\end{lemma}
\begin{main}\label{finitetype} Let $\Gamma$ be an $N=1$ ADE quiver. Let
  $\mathcal{B}_\Gamma= \{r_i(P'_j(x))\},$ where $r_i\in {\mathfrak
    W}_\Gamma$ and $P'_j,\, j\in V_\Gamma$ are the polynomials defined
  in relation (\ref{relation}).  Assume no element in
  $\mathcal{B}_\Gamma$ has a multiple root. If $(*)$ holds, then $N=1$
  ADE quivers have finite representation type.
\end{main}
\begin{proof} Let $V$ be a simple representation of an $N=1$ ADE quiver.
Let $\mathcal A=\{i|V(i)\not=0\}.$  We can assume that $\mathcal A$ is
connected.   Otherwise, $V$ would be decomposable.    We apply the
forgetful functors to $V$ to get an (+)-accessible (resp.
(-)-accessible) diagram (no loop)

\[\xymatrix
{V(1)&\cdots\ar[l]& V(k)\ar[l]&\cdots\ar[l]& V(n)\ar[l]\\
&& V(l)\ar[u]\\
}
\]
\[\xymatrix
{V(1)\ar[r]& \cdots\ar[r]& V(k)\ar[r]\ar[d]&\cdots\ar[r]& V(n)\\
&& V(l)\\
}
\]
(For the type A case, $V(l)=0.$)

Let $c=\sigma_n\cdots \sigma_2\sigma_1.$  By \cite{BGP}, there exists $k$
such that $c^k(\dim V)\ngtr 0.$  By $(*)$, Lemma~\ref{poly1} and
Lemma~\ref{poly2},  we know that $\sum_{i} \dim V(i).P'_i(x)$ is the only
element in ${\mathcal A}_\Gamma$  which vanishes at $\lambda.$   By
Corollary~\ref{cor} and Proposition~\ref{newrelation},  this implies that
there exist $\beta_1,\cdots, \beta_l$ and a simple representation
$L_{\beta_{l+1}}$ which is concentrated at a vertex of $\Gamma$  such that
\[V=F_{\beta_1}\cdots F_{\beta_{k}}(L_{\beta_{k+1}})\]
$V$ corresponds to the  positive root
\[\dim V=\sigma_{\beta_1}\cdots \sigma_{\beta_k}(\overline{\beta_{k+1}})\]
where
\[
(\overline{\beta_{k+1}})=\begin{cases} 0 \quad \mbox{if} \quad i\not=k+1\\
                                1  \quad \mbox{if} \quad i=k+1\\
                  \end{cases}
\]
Since  the usual ADE quiver only has finitely many positive roots, $N=1$
ADE quivers have finite representation type. This finishes the proof of
the theorem.
\end{proof}
From the above Main~Theorem, one can get the following
Proposition~\ref{finite-1corr}.
\begin{corollary}\label{finite-1corr} Let $\Gamma$ be an $N=1$ ADE
quiver.  Let $\mathcal {B}_\Gamma=\{r_i  (P'_j(x))  | r_i  \in \mathfrak
{W}_\Gamma\}$, where $\mathfrak {W}_\Gamma$ is the Weyl group of $\Gamma$
and $P'_j$ is the polynomial defined on relation~(\ref{relation}).  Assume  
each element in $\mathcal
{B}_\Gamma$ has simple roots. If $(*)$ holds,  then there is a
finite-to-one correspondence between simple representations of $N=1$ ADE
quivers and the positive roots of ADE Dynkin diagram.
\end{corollary}
\begin{proof} We know that $\mathcal{B}_\Gamma$ has only finitely many
elements.  Each element of $\mathcal{B}_\Gamma$ which is in fact a
polynomial has only finitely many simple roots. By our Main Theorem, each
root of an element in $\mathcal{B}_\Gamma$ corresponds with a simple
representation.  Hence, the desired result follows.
\end{proof}

\section{A correspondence between indecomposable representations and ADE configuration of curves.}\label{section5}

This section will give a correspondence between indecomposable
representations of an N=1 ADE quiver and ADE configuration of curves.

To make the presentation easier to follow,
we state the following Lemma from \cite[(1.1),(1.14)]{Re83} here.

{\bf Reid's Lemma.}
 {\em Let $\pi\colon Y \to X$ be a resolution of an isolated Gorenstein
  threefold singularity $P \in X$.  Suppose that the exceptional set
  of $\pi$ has pure dimension 1.  Let $X_0$ be a generic hyperplane
  section of $X$ which passes through $P$.  Then $X_0$ has a rational
  double point at $P$.
  
  Moreover, if $X_0$ is any hyperplane section through $P$ with a
  rational double point, and $Y_0$ is its proper transform, then $Y_0$
  is normal, and the minimal resolution $Z_0 \to X_0$ factors through
  the induced map $\pi |_{Y_0}\colon Y_0 \to X_0$.}

The exceptional curve of the minimal resolution $Z_{0}\to X_{0}$ is a
rational curve whose dual graph is $\Gamma $. Thus the components of the
exceptional divisor are rational curves $C_{e_{i}}$ which are in one to one
correspondence with the simple roots $e_{1},\dots ,e_{n}$ ($n$ is the rank
of the root system).

Let $X$ be an ADE fibration with base $\mathbb C,$ and $Y$ be the
small resolution of $X.$ Let $\pi: Y\to X$ be the blow up map. An
``ADE configuration of curves'' in $Y$ is a 1 dimensional connected
projective scheme $C\subset Y,$ such that

\begin{enumerate}
\item $\exists \bar{S}\subset Y,\, C\subset \bar{S}$
\item letting $S=\pi(\bar S),$ then $\bar{S}\to S$ is a resolution of
  $ADE$ singularities with exceptional scheme $C.$
\end{enumerate}

We need the following proposition which is essentially part 3 of Theorem 1 in
\cite{KM92}.
\begin{proposition}\label{prop: components of disc locus of Res}
The irreducible components of the discriminant divisor $\mathfrak{D} \subset
\operatorname{Res} (\Gamma )$ are in one to one correspondence with the
positive roots of $\Gamma $. Under the identification of
$\operatorname{Res} (\Gamma )$ with the complex root space $U$, the
component $\mathfrak{D} _{v}$ corresponding to the positive root $v=\sum
_{i=1}^{n}a_{i}e_{i}$ is $v^{\perp
}\subset U $,  ie  the hyperplane perpendicular to $v$.

Moreover, $\mathfrak{D} _{v}$ corresponds exactly to those deformations of
$Z_{0}$ in $\mathcal{Z}$ to which the curve 
\[
C_{v}:=\bigcup
_{i=1}^{n}a_{i}C_{e_{i}}
\]
lifts. For a generic point $t\in \mathfrak{D} _{v}$,
the corresponding surface $\mathcal{Z}_{t}$ has a single smooth $-2$ curve
in the class $\sum _{i=1}^{n}a_{i}[C_{e_{i}}]$ thus there is a small
neighborhood $B$ of $t$ such that the restriction of $\mathcal{Z}$ to $B$
is isomorphic to a product of $\C ^{n-1}$ with the semi-universal
family over $\operatorname{Res} (A_{1})$.
\end{proposition}
Let $X$ be an ADE fibration with base $\mathbb C,$ and $Y$ be the
  small resolution of $X.$
  
 The following example gives a concrete correspondence between
  the ADE configuration of curves in $Y$ and the indecomposable
  representations of an N=1 ADE quiver in the $A_2$ case.

\begin{example} Let $X$ be defined by 
\[
A_2: \, xy+(z+t_1(t))(z+t_2(t))(z+t_3(t))=0
\]
with 
\[
t_1(t)+t_2(t)+t_3(t)=0.
\]
In the following Table~\ref{table1},
``Curve'' means an ADE configuration of curves which is in the
exceptional set of the fibration. I use ``$\dim \,V$'' to denote the
dimension vector of an indecomposable representation of the N=1 ADE quiver.

Table~\ref{table1} can be explained in the following way.  If
$t_1(\lambda)=t_2(\lambda)\ne t_3(\lambda)$ for some $\lambda,$  then by
Proposition~\ref{prop: components of disc locus of Res}, there exists an ADE
configuration of curves $C\subset Y.$ By \cite{CKV}, we know that
$P'_1(\lambda)=t_1(\lambda)-t_2(\lambda)=0.$ Hence, by our Main
Theorem, there exists an indecomposable representation $V$ of an N=1
ADE quiver with $\dim V=(1,0)$ that corresponds to $C.$

The other cases are similar, we omit them.

\begin{table}
\begin{center}
\begin{tabular}{|c|c|c|c|c|}
\hline 
Condition &  Singularity & Curve & $\dim\, V$ \\
\hline
$t_1(\lambda)=t_2(\lambda)\ne t_3(\lambda)$ & $A_1$ & $\xymatrix{   &\\
\ar@{-}[ur]^{{\mathbb P}^1}\\
}$ & $\xymatrix@R=0pt{
1 & 0\\
\bullet\ar@{-}[r] &\bullet\\}$\\
\hline
$t_1(\lambda)\ne t_2(\lambda)=t_3(\lambda)$ & $A_1$  & $\xymatrix{   &\\
\ar@{-}[ur]^{{\mathbb P}^1}\\}$ & $\xymatrix@R=0pt{
0 & 1\\
\bullet\ar@{-}[r] &\bullet\\}$ \\
\hline
$t_1(\lambda)=t_3(\lambda)\ne t_2(\lambda)$  & $A_2$  & $\xymatrix{
        &      &            &             &        \\
              &    &             &        \\
\ar@{-}[uurr]^{{\mathbb P}^1} &              &\\
&  & \ar@{-}[uull]^{{\mathbb P}^1}\\ }$ & $\xymatrix@R=0pt{
1 & 1\\
\bullet\ar@{-}[r] &\bullet\\}$\\
\hline
$t_1(\lambda)\ne t_2(\lambda)\ne t_3(\lambda)$ & -  & - &-\\
\hline
\end{tabular}
\end{center}

\caption{}
\label{table1}
\end{table}
\end{example}

This example will be generalized in the following theorem.

\begin{theorem}\label{representations and curves}
  Let $X$ be a ADE fibration corresponding to $\Gamma,$ with base $\mathbb C.$  Let $Y$ be
 a small resolution of $X.$   Let
  $\mathcal {B}_\Gamma=\{r_i (P'_j(x)) | r_i \in \mathfrak
  {W}_\Gamma\}$, where $\mathfrak {W}_\Gamma$ is the Weyl group of
  $\Gamma$ and $P'_j$ is the polynomial defined in
  relation~(\ref{relation}).  Assume no element in $\mathcal
  {B}_\Gamma$ has multiple roots and assume $(*)$ holds. Then there exists
  a 1-1 correspondence between the indecomposable representations of the
  $N=1$ ADE quiver and the ADE configuration of curves in $Y$.
\end{theorem}
\begin{proof}
  By Pinkham \cite{Pinkham} and Katz-Morrison \cite{KM92}, we have the
  following commutative diagram
\[
\xymatrix{
Y\ar[r]\ar[d]^{\pi} & \mathcal{Y}\ar[d]\\
X\ar[r]\ar[d]^{\varphi} & \mathcal{X}\ar[d]\\
\C\ar[r]^f &       \C^n\\
}
\]
where $\C$ denotes the set of complex numbers and $\mathcal Y$ denotes
the $\C^*$-equivariant simultaneous resolution $\mathcal{Y}\to
\mathcal{X}$ inducing $Y_0\to X_0.$ For an indecomposable
representation $V$ of the $N=1$ ADE quiver $\Gamma,$  we have
\begin{equation}\label{equation_a}
\sum \dim V(i)\cdot P'_i(\lambda)=0
\end{equation}
for some $\lambda.$ The dimension vector $\left(\dim\,
  V(i)\right)_{i\in V_\Gamma}$ will correspond to a positive root
$\rho.$ By (\ref{A_n}), (\ref{D_n}), and (\ref{E_n}), we can express
$P'_i(x),\, i=1,\cdots , n$ in terms of $t_i,\, i=1,\cdots , n.$ By
Proposition ~\ref{prop: components of disc locus of Res} or part 3 of
Theorem 1 in Katz-Morrison \cite[pp.~467]{KM92}, (\ref{equation_a})
will give an equation for $\rho^\perp.$ Hence
$f(\lambda)=\left(t_i(\lambda)\right)_{i\in V_\Gamma} \in \rho^\perp.$
It follows from Proposition ~\ref{prop: components of disc locus of
  Res} that there exists an ADE configuration of curves $C_\rho\subset
\pi^{-1}(\lambda)\subset Y.$

Conversely, for an ADE configuration of curves $C\subset Y,$ we have
that $\varphi\circ \pi(C)=\lambda \in \C$ ( Since $\pi$ is projective,
$\varphi\circ \pi(C)$ is projective in $\C.$ It follows that
$\varphi\circ\pi(C)$ is a finite subset of $\C.$ Since $C$ is
connected, $\varphi\circ \pi(C)$ is connected in $\C.$ Hence
$\varphi\circ \pi$ is a point in $\C.$ ) Moreover, $\pi(C)$ is a point
in $X$ (By Katz-Morrison \cite{KM92}, we know that $\mathcal{X}$ is
affine.  Hence $\pi(C)$ is a point in $X.$) By Proposition ~\ref{prop:
  components of disc locus of Res}, we know that $f(\lambda)\in
\rho^\perp$ for some positive root $\rho. $ Since we assume that each
element in ${\mathcal B}_\Gamma$ has simple roots and $(*)$ holds, $C$
corresponds to a unique positive root $\rho$.  We can express $\rho$
as $\rho=\sum a_i\cdot \rho_i$ where $\rho_i$ is a simple positive
root.  From our Main Theorem, we can construct a simple representation
$V$ of $N=1$ ADE quiver $\Gamma$ which corresponds to the positive
root $\rho$ by applying the reflection functors.  This finishes the
proof of Theorem \ref{representations and curves}.

\end{proof}


\end{document}